\documentclass{article}
\usepackage{graphicx}
\usepackage{amsmath}
\usepackage{amsfonts}
\usepackage{amssymb}

\newtheorem{theorem}{Theorem}

\newtheorem{corollary}{Corollary}

\newtheorem{definition}{Definition}

\newtheorem{lemma}{Lemma}

\newtheorem{proposition}{Proposition}
\newtheorem{remark}{Remark}

\begin{document}

\title{NON EXISTENCE\ OF\ COMPLEX\ STRUCTURES\\
\ ON FILIFORM LIE\ ALGEBRAS\ }

\author{Michel Goze \thanks{%
corresponding author: e-mail: M.Goze@uha.fr} \and Elisabeth
Remm\\
Universit\'{e} de Haute Alsace, F.S.T.\\
4, rue des Fr\`{e}res Lumi\`{e}re - 68093 MULHOUSE - France}
\date{}
\maketitle

\begin{abstract}
The aim of this work is to prove the nonexistence of complex structures over
nilpotent Lie algebras of maximal class (also called filiform).
\end{abstract}

\section{Preliminaries}

\bigskip The study of invariant complex structures on real connected Lie
Groups is reduced to the study of linear operators $J$ on the corresponding
Lie algebra which satisfies the Nijenhuis condition and $J^{2}=-Id$.\ When
the Lie algebra is even dimensional, real and reductive, the existence of
such structures follows from the work of Morimoto [1]. On the other  hand,
there exist solvable and nilpotent Lie algebras which are not provided  with
a complex structure. Some recent works present classifications of  nilpotent
or solvable Lie algebras with complex structures in small  dimension
(dimension four for the solvable case [2] and six for the nilpotent case
[3]).\ As every six dimensional nilpotent Lie algebra does not admit
complex structure, we are conducted to determine the classes of nilpotent
Lie algebra which are not provided with such a structure.

Let $G$ be a real Lie group and $\frak{g}$ its Lie algebra.

\begin{definition}
An invariant complex structure on $G$ is an endomorphism $J$ of $\frak{g}$
such that

(1) $J^{2}=-Id$

(2) $\left[ JX,JY\right] =\left[ X,Y\right] +J\left[ JX,Y\right] +J\left[
X,J(Y)\right] ,$ $\forall X,Y\in \frak{g}$
\end{definition}

(The second condition is called the Nijenhuis condition of integrability).

For simplify, we will say that $J$ is an invariant complex structure over
the  Lie algebra $\frak{g}$.

The underlying real vector space $V$ to the Lie algebra $\frak{g}$ can be
provided with a complex vector space structure by putting 
\[
\left( a+ib\right).v=a.v+b.J(v) 
\]
$\forall v \in V, \forall a,b \in \Bbb{R}.$ We note by $V_{J}$ this complex
vector space.\ We have 
\[
\dim _{\Bbb{C}}V_{J}=\frac{1}{2}\dim _{\Bbb{R}}V=\frac{1}{2}\dim _{\Bbb{R}}%
\frak{g}.
\]

\begin{definition}
The complex structure $J$ is called bi-invariant if it satisfies

(3) $\left[ J,adX\right] =0$ , $\forall X\in \frak{g}$.
\end{definition}

Let us remark that (3) implies (2).

If $J$ is bi-invariant, then $V_{J}$ is a complex Lie algebra, noted $\frak{g%
}_{J}$, because in this case we have 
\[
\left[ \left( a+ib\right) X,\left( c+id\right) Y\right] =\left( a+ib\right)
\left( c+id\right) \left[ X,Y\right] 
\]

\section{Decomposition associated to a complex structure}

Let $J$ be an invariant complex structure on the real Lie algebra $\frak{g}$%
.\ We can extend the endomorphism $J$ in a natural way on the
complexification $\frak{g}_{\Bbb{C}}=\frak{g}{\otimes }_{\Bbb{R}}\Bbb{C}$ of 
$\frak{g}$.\ It induces on $\frak{g}_{\Bbb{C}}$ a direct vectorial sum 
\[
\frak{g}_{\Bbb{C}}=\frak{g}_{\Bbb{C}}^{i}\oplus \frak{g}_{\Bbb{C}}^{-i} 
\]
where 
\[
\frak{g}_{\Bbb{C}}^{\varepsilon i}=\left\{ X\in \frak{g}_{\Bbb{C}}\quad
/\quad J(X)=\varepsilon iX\right\} \quad , \quad \varepsilon =\pm 1 
\]

The Nijenhuis condition (2) implies that $\frak{g}_{\Bbb{C}}^{\varepsilon i}$
is a complex subalgebra of $\frak{g}_{\Bbb{C}}^{{}}$.\ If $\sigma$ denotes
the conjugation on $\frak{g}_{\Bbb{C}}$ given by $\sigma \left( X+iY\right)
=X-iY$, then 
\[
\frak{g}_{\Bbb{C}}^{-i}=\sigma (\frak{g}_{\Bbb{C}}^{i}). 
\]

\begin{proposition}
Let $\frak{g}$ be a $2n$-dimensional real Lie algebra.\ It is provided with
an invariant complex structure if and only if the complexification $\frak{g}%
_{\Bbb{C}}$ admits a decomposition

\[
\frak{g}_{\Bbb{C}}=\frak{h}\oplus \sigma (\frak{h})
\]
where $\frak{h}$ is a $n$-dimensional complex subalgebra of $\frak{g}_{\Bbb{C%
}}$.
\end{proposition}

If $J$ is bi-invariant, condition (3) implies that $\frak{g}_{\Bbb{C}%
}^{\varepsilon i}$ is an ideal of $\frak{g}_{\Bbb{C}}$.\ In fact, if $X\in 
\frak{g}_{\Bbb{C}}^{\varepsilon i}$ and $Y\in \frak{g}_{\Bbb{C}%
}^{-\varepsilon i}$ we have 
\[
J\left[ X,Y\right] =\left[ JX,Y\right] =\varepsilon i\left[ X,Y\right]
=\left[ X,JY\right] =-\varepsilon i\left[ X,Y\right] 
\]
Then 
\[
\left[ X,Y\right] =0. 
\]

\begin{proposition}
Let $\frak{g}$ be an $2n$-dimensional real Lie algebra.\ Then $\frak{g}$
admits a bi-invariant complex structure if and only if the complexification $%
\frak{g}_{\Bbb{C}}$ is a direct sum of ideals $I$ and $\sigma \left(
I\right) $ : 
\[
\frak{g}_{\Bbb{C}}=I\oplus \sigma (I)
\]
\end{proposition}

Of course every $n$-dimensional complex Lie algebra $\frak{h}$ comes from an 
$2n$-dimensional real Lie algebra endowed with a bi-invariant complex
structure.

\section{Bi-invariant complex structures}

\subsection{Nilpotent case}

Let $\frak{g}$ be a $2n$-dimensional real nilpotent Lie algebra with a
bi-invariant complex structure. Then $\frak{g}_{\Bbb{C}}=I\oplus \sigma (I)$
where $I$ is a $n$-dimensional complex ideal of $\frak{g}_{\Bbb{C}}$.\ This
describes entirely the structure of the complexifications $\frak{g}_{\Bbb{C}%
} $ of the real Lie algebras $\frak{g}$ provided with bi-invariant complex
structure.

Let $c(\frak{n})$ the characteristic sequence of the complex nilpotent
algebra $\frak{n}$ ([4]).\ It is defined by 
\[
\begin{array}{ll}
c(\frak{n})= & Max\{c(X)\quad \mid \quad X\in \frak{n}-\left[ \frak{n},\frak{%
n}\right] \}
\end{array}
\]
where $c(X)=\left( c_{1}(X),...,c_{k}(X),1\right) $ is the sequence of
similitude invariants of the nilpotent operator $adX$.\ We deduce that 
\[
c\left( \frak{g}_{\Bbb{C}}\right) =\left(
c_{1},c_{1},c_{2},c_{2},...,1,1\right) 
\]
As the Jordan normal form of the nilpotent operators $adX$ does not depend
of the  field of scalars, we have

\begin{proposition}
If $\frak{g}$ is an $2n$-dimensional real nilpotent Lie algebra which admits
a bi-invariant complex structure, then its characteristic sequence is of
type 
\[
\left( c_{1},c_{1},c_{2},c_{2},...,c_{j},c_{j},...,1,1\right) .
\]
\end{proposition}

A nilpotent Lie algebra is called of maximal class (or filiform in  Vergne's
terminology [5],[6]) if the descending sequence of derived ideals $\mathcal{C}%
^{i}\frak{g}$ satisfies : 
\[
\begin{array}{l}
\dim \mathcal{C}^{1}\frak{g}=\dim \frak{g}-2 \\ 
\dim \mathcal{C}^{i}\frak{g}=\dim \frak{g}-i-1\quad \forall i\geq 2
\end{array}
\]

\begin{corollary}
There is no bi-invariant complex structure over a filiform Lie algebra.
\end{corollary}

In fact, the characteristic sequence of a filiform Lie algebra is $(2n-1,1).$
Thus, by proposition 3, a filiform Lie algebra cannot admit a bi-invariant
complex structure.

Let be $\left\{ U_{j}\right\} $ and $\left\{ V_{j}\right\} $ the basis of
the complex ideals $I$ and $\sigma (I)$. Respect to $\frak{g}$, we can write 
$U_{l}=X_{l}-iY_{l}$, $V_{l}=X_{l}+iY_{l}$, where $\left\{
X_{1},...,X_{n},Y,...,Y_{n}\right\} $ is a real basis of $\frak{g}$.\ As $%
\left[ U_{l},V_{l}\right] =0$, we obtain 
\[
\left[ X_{l}-iY_{l},X_{k}+iY_{k}\right] =\left[ X_{l},X_{k}\right] +\left[
Y_{l},Y_{k}\right] +i\left( \left[ X_{l},Y_{k}\right] -\left[
Y_{l},X_{k}\right] \right) =0 
\]
Thus 
\[
\begin{array}{ll}
\left[ X_{l},X_{k}\right] =-\left[ Y_{l},Y_{k}\right] &  \\ 
\left[ X_{l},Y_{k}\right] =\left[ Y_{l},X_{k}\right] & k=1,...,n\quad ;
\quad l=1,...,n.
\end{array}
\]
Likewise $\left[ U_{l},U_{k}\right] \in I$ and $\left[ V_{l},V_{k}\right]
\in \sigma \left( I\right) $ imply 
\[
\begin{array}{l}
\left[ X_{l},X_{k}\right] -\left[ Y_{l},Y_{k}\right] \in \Bbb{C}\left\{
X_{1},...,X_{k}\right\} \\ 
\left[ X_{l},Y_{k}\right] +\left[ Y_{l},Y_{k}\right] \in \Bbb{C}\left\{
Y_{1},...,Y_{k}\right\}
\end{array}
\]
Suppose that $I$ admits a real basis.\ In this case we have 
\[
\begin{array}{lll}
\left[ X_{l},X_{k}\right] ={\sum }_{j=1}^{n}C_{lk}^{j}X_{j} & , & \left[
Y_{l},Y_{k}\right] =-{\sum }_{j=1}^{n}C_{lk}^{j}Y_{j} \\ 
\left[ X_{l},Y_{k}\right] ={\sum }_{j=1}^{n} D_{lk}^{j}Y_{j} & , & \left[
Y_{l},X_{k}\right] ={\sum }_{j=1}^{n}D_{lk}^{j}Y_{j}
\end{array}
\]
with $C_{lk}^{j}$ and $D_{lk}^{j}$ in $\Bbb{R}$. Let be $\frak{g}_{0}=\Bbb{R}%
\left\{ X_{1},...,X_{n}\right\} $ and $\frak{g}_{1}=\Bbb{R}\left\{
Y_{1},...,Y_{n}\right\} $.\ Then the previous relations show that $\frak{g}=%
\frak{g}_{0}\oplus \frak{g}_{1}$, and 
\[
\left\{ 
\begin{array}{l}
\left[ \frak{g}_{0},\frak{g}_{0}\right] \subset \frak{g}_{0} \\ 
\left[ \frak{g}_{1},\frak{g}_{1}\right] \subset \frak{g}_{0} \\ 
\left[ \frak{g}_{0},\frak{g}_{1}\right] \subset \frak{g}_{1}
\end{array}
\right. 
\]
which gives a structure of $\Bbb{Z}_{2}$-graded Lie algebra on $\frak{g}.$

\subsection{On the classification of nilpotent Lie algebras with
bi-invariant complex structures}

If $\frak{g}$ is a $2n$-dimensional real nilpotent Lie algebra with a
bi-invariant complex structure, then $\frak{g}_{\Bbb{C}}=I\oplus \sigma (I)$%
.\ Then the classification of associated complexification $\frak{g}_{\Bbb{C}}
$ corresponds to the classification of complex nilpotent Lie algebra $I$.\
By this time, this classification is known only up to dimension 7, and for
special classes up to dimension 8. To obtain a general list of these
algebras is therefore a hopeless pretention.\ From [3] we can extract the
list for the 6-dimensional case. From [4] we can present the
classification of complexifications $\frak{g}_{\Bbb{C}}$ of real nilpotent
Lie algebras with bi-invariant complex structure up the dimension 14. Here
we will present briefly the classification in the real case up dimension 8.

\textbf{Dimension 2} : $\frak{g}_{2}^{1}$ is abelian.

\textbf{Dimension 4} : $\frak{g}_{4}^{1}$ is abelian.

In fact $\frak{g}_{\Bbb{C}}=I\oplus \sigma (I)$ and $I$ is an abelian ideal.

\textbf{Dimension 6} : - $\frak{g}_{6}^{1}$ is abelian.

\textbf{- }$\frak{g}_{6}^{2}$ : $\left\{ 
\begin{array}{l}
\left[ X_{2},X_{3}\right] =-\left[ Y_{2},Y_{3}\right] =X_{1} \\ 
\left[ X_{2},Y_{3}\right] =\left[ Y_{2},X_{3}\right] =Y_{1}
\end{array}
\right. $

\begin{remark}
1.\ The complexification of $\frak{g}_{6}^{2}$ is isomorphic to $\frak{h}%
_{3}\oplus \frak{h}_{3}$ where $\frak{h}_{3}$ is the $3$-dimensional
Heisenberg Lie algebra.

2. Let us note that the real nilpotent Lie algebra $\frak{h}_{3}\oplus \frak{%
h}_{3}$ doesn't have a bi-invariant complex structure.\ In fact if $\left\{
X_{1},X_{2},X_{3},X_{4},X_{5},X_{6}\right\} $ is a basis of $\frak{h}%
_{3}\oplus \frak{h}_{3}$ with brackets 
\[
\left[ X_{1},X_{2}\right] =X_{3},\quad \left[ X_{4},X_{5}\right] =X_{6},
\]
then every isomorphism $J$ which commutes with $adX$ for every $X$ satisfies 
$J\left( X_{3}\right) =\alpha X_{3}$ and $J\left( X_{6}\right) =\beta X_{6}$%
. As $J^{2}=-Id$, we will have $\alpha ^{2}=\beta ^{2}=-1$ and $\alpha
,\beta \in \Bbb{R}$.

3. The Lie algebra $\frak{g}_{6}^{2}$ is $2$-step nilpotent and it is a Lie
algebra of type $H$ (see [7]).
\end{remark}

\textbf{Dimension 8} : - $\frak{g}_{8}^{1}$ is abelian.

\textbf{- }$\frak{g}_{8}^{2}=\frak{g}_{6}^{2}\oplus \frak{g}_{2}^{1}$

- $\frak{g}_{8}^{3}:\left\{ 
\begin{array}{lll}
\left[ X_{2},X_{4}\right] =X_{1} & ; & \left[ Y_{2},Y_{4}\right] =-X_{1} \\ 
\left[ X_{2},Y_{4}\right] =+Y_{1} & ; & \left[ Y_{2},X_{4}\right] =+Y_{1} \\ 
\left[ X_{3},X_{4}\right] =X_{2} & ; & \left[ Y_{3},Y_{4}\right] =-X_{2} \\ 
\left[ X_{3},Y_{4}\right] =Y_{2} & ; & \left[ Y_{3},X_{4}\right] =Y_{2}
\end{array}
\right. $

- $\frak{g}_{8}^{4}:\left\{ 
\begin{array}{lll}
\left[ X_{2},X_{3}\right] =X_{1} & ; & \left[ Y_{2},Y_{3}\right] =-X_{1} \\ 
\left[ X_{2},Y_{3}\right] =Y_{1} & ; & \left[ Y_{2},X_{3}\right] =Y_{1} \\ 
\left[ X_{2},X_{4}\right] =Y_{1} & ; & \left[ Y_{2},Y_{4}\right] =-Y_{1} \\ 
\left[ X_{2},Y_{4}\right] =-X_{1} & ; & \left[ Y_{2},X_{4}\right] =-X_{1}
\end{array}
\right. $

\subsection{On the classification of solvable Lie algebras with bi-invariant
complex structures}

As known, there is, up an isomorphism, only one $2$-dimensional non
nilpotent solvable Lie algebra, denoted by $\frak{r}_{2}^{2}$\ and defined
by $\left[ X_{1},X_{2}\right] =X_{1}$ over the basis $\left\{X_{1},X_{2}%
\right\} $ .  This Lie algebra does not admit a bi-invariant complex
structure.\ Thus, a non abelian solvable Lie algebra admitting it is at
least four dimensional.

If we want classify all the $4$-dimensional solvable Lie algebras admitting 
bi-invariant complex structures we can use the Dozias'classification, which
can be found in [8] (see also [9]).

From this classification, we can affirm :

\begin{theorem}
: Every $4$-dimensional real solvable Lie algebras which admits a
bi-invariant complex structure is isomorphic to one of the following :

- $\frak{r}_{4}^{1}=\frak{g}_{4}^{1}$ the abelian Lie algebra.

- $\frak{r}_{4}^{2}=\left\{ 
\begin{array}{c}
\lbrack X_{1},X_{3}]=X_{3} \\ 
\lbrack X_{1},X_{4}]=X_{4} \\ 
\lbrack X_{2},X_{3}]=-X_{4} \\ 
\lbrack X_{2},X_{4}]=X_{3}
\end{array}
\right. $
\end{theorem}

In the last case, any bi-invariant complexe structures satisfy 
\begin{eqnarray*}
J(X_{1}) &=&aX_{1}+bX_{2}, \\
J(X_{2}) &=&-bX_{1}+aX_{2} \\
J(X_{3}) &=&aX_{3}+bX_{4}, \\
J(X_{4}) &=&-bX_{3}+aX_{4}.
\end{eqnarray*}
Its complexification $\frak{r}_{4\Bbb{C}}^{2}$ is isomorphic to $\frak{r}%
_{2}^{2}\times \frak{r}_{2}^{2}.$ Let us note that the real Lie algebra $%
\frak{r}_{2}^{2}\times \frak{r}_{2}^{2}$ which has the same complexification
as $\frak{r}_{4}^{2}$ is not provided with a bi-invariant complex structure.

\section{Non Existence of invariant complex structures over nilpotent Lie
algebras of maximal class}

In section 3, we have given the definition of nilpotent Lie algebras of
maximal class (or filiform Lie algebra).The simplest example is given by the
following $n$-dimensional Lie algebra, denoted $L_{n}$ : 
\[
\left[ X_{1},X_{i}\right] =X_{i+1},\quad i=2,...,n-1
\]
where the nondefined brackets are zero or obtained by antisymmetry.

Let $\frak{n}$ be a $n$-dimensional filiform Lie algebra. Then there exists a
basis $\{X_{1},X_{2},...,X_{n}\}$ which is adapted to the flag 
\[
\frak{g}\supset \mathcal{C}^{1}\frak{g}\supset \mathcal{C}^{2}\frak{g}%
\supset ...\supset \mathcal{C}^{n-1}\frak{g}=\left\{ 0\right\} 
\]
with $\dim \mathcal{C}^{1}\frak{g}=n-2$ , $\dim \frac{\mathcal{C}^{i}\frak{g}%
}{\mathcal{C}^{i-1}\frak{g}}$ =1, $i\geq 1$, and which satisfies : 
\[
\left\{ 
\begin{array}{l}
\left[ X_{1},X_{i}\right] =X_{i+1},\quad i=2,...,n-1 \\ 
\left[ X_{i},X_{j}\right] ={\sum }_{k \geq i+j}C_{ij}^{k}X_{k}
\end{array}
\right. 
\]
The change of basis $Y_{1}=X_{1}$, $Y_{i}=tX_{i}$, $i\geq 2$, $t\neq 0$,
shows that this Lie algebra is isomorphic to the following : 

\[
\left\{ 
\begin{array}{l}
\left[ X_{1},X_{i}\right] =X_{i+1},\quad i=2,...,n-1 \\ 
\left[ X_{i},X_{j}\right] =t{\sum }C_{ij}^{k}X_{k}
\end{array}
\right. 
\]

\subsection{Invariant complex structures and filiform Lie algebras}

\begin{proposition}
The real nilpotent filiform Lie algebra $L_{2n}$ does not admit an invariant
complex structure.
\end{proposition}

{\bf Proof }. Let $T$ be a linear isomorphism of the real vector space $L_{2n}$
satisfying the Nijenhuis condition : 
\[
\left[ T(X),T(Y)\right] =\left[ X,Y\right] +T\left[ X,T(Y)\right] +T\left[
T(X),Y\right] 
\]
where $\left[ ,\right] $ is the bracket of $L_{2n}$. Consider the basis $%
\left\{ X_{1},...,X_{2n}\right\} $ of $L_{2n}$ satisfying : 
\[
\left\{ 
\begin{array}{ll}
\left[ X_{1},X_{i}\right] =X_{i+1}, & i=2,...,2n-1 \\ 
\left[ X_{i},X_{j}\right] =0, & i,j\neq 1.
\end{array}
\right. 
\]
We have 
\begin{eqnarray*}
\left[ T(X_{2n-1}),T(X_{2n})\right]  &=&\left[ X_{2n-1},X_{2n}\right]
+T\left[ X_{2n-1},T(X_{2n})\right] +T\left[ T(X_{2n-1}),X_{2n}\right]  \\
&=&T\left[ X_{2n-1},T(X_{2n})\right] 
\end{eqnarray*}
As 
\[
\left[ X_{2n-1},X_{1}\right] =-X_{2n}
\]
we obtain 
\[
\left[ T(X_{2n-1}),T(X_{2n})\right] =T\left[ X_{2n-1},T(X_{2n})\right]
=-aT(X_{2n})
\]
where 
\[
T(X_{2n})=aX_{1}+\sum_{i\geq 2}a_{i}X_{i}.
\]
The nilpotency of $L_{2n}$ implies that the constant $a$ which appears as an
eigenvalue of $ad$($T(X_{2n-1}))$ is zero.\ Then 
\[
\left[ T(X_{2n-1}),T(X_{2n})\right] =0.
\]
This implies that 
\[
\left[ T(X_{i}),T(X_{2n})\right] =T\left[ X_{i},T(X_{2n})\right] =0
\]
for $i=2,...,2n.$ If $T(X_{2n})\notin Z(L_{2n})=\Bbb{R}\{X_{2n}\},$ then $%
T(X_{i})=\sum_{j\geq 2}a_{ij}X_{j}$ for $j=2,...,2n.$ As $T$ is non
singular, necessarily we have 
\[
T(X_{1})=a_{11}X_{1}+\sum_{i\geq 2}a_{i1}X_{i}
\]
with $a_{11}\neq 0.$ The condition $T^{2}=-Id$ implies $a_{11}^{2}=-1.$ As $%
a_{11}^{{}}\in \Bbb{R}$, this constitutes a contradiction. Thus $%
T(X_{2n})\in Z(L_{2n})=\Bbb{R}\{X_{2n}\},$ and it follows that $%
T(X_{2n})=\alpha X_{2n}.$ As above, we can conclude $\alpha ^{2}=-1.$ This
case is also excluded and the proposition is proved.

We will show that the non existence of invariant complexe structure on the
model filiform $L_{2n}$ is always true in any deformation of this algebra.

\begin{theorem}
There does not exist invariant complex structure over a real filiform Lie
algebra.
\end{theorem}

{\bf Proof}. Let $\frak{g}$ be a $2n$-dimensional real filiform Lie algebra and
let $\frak{g}_{\Bbb{C}}$ be its complexification. If there a complex
structure $J$ on $\frak{g}$, then $\frak{g}_{\Bbb{C}}$ admits the following
decomposition 
\[
\frak{g}_{\Bbb{C}}=\frak{g}_{1}\oplus \sigma (\frak{g}_{1}) 
\]
where $\frak{g}_{1}$ is a $n$-dimensional sub-algebra of$\frak{\ g}_{\Bbb{C}%
} $ and $\sigma $ the conjugated automorphism. As the Lie algebra $\frak{g}_{%
\Bbb{C}}$ is filiform, there is an adapted basis $\{X_{1},..,X_{2n}\}$
satisfying 
\[
(*)\left\{ 
\begin{array}{l}
\lbrack X_{1},X_{i}]=X_{i+1},\quad 2\leq i\leq 2n-1 \\ 
\lbrack X_{2},X_{3}]=\sum_{k\geq 5}C_{23}^{k}X_{k} \\
{\cal{C}} ^i(\frak g)=\Bbb{R}\{X_{i+2},...,X_n\}
\end{array}
\right. 
\]
In particular we have 
\[
\dim \frac{\frak{g}_{\Bbb{C}}}{\mathcal{C}^{1}(\frak{g}_{\Bbb{C}})}=2,\quad
\dim \frac{\mathcal{C}^{i}(\frak{g}_{\Bbb{C}})}{\mathcal{C}^{i+1}(\frak{g}_{%
\Bbb{C}})}=1\quad , \quad i\geq 1 
\]
The ordered sequence of the dimension of Jordan blocks of the nilpotent
operator $adX_{1}$ is $(n-1,1)$. Such a vector is called characteristic
vector.

\begin{lemma}
Every characteristic vector can be written as $Y=\alpha X_{1}+U$ where $U$
is in the complex vector space generated by $\{X_{2},...,X_{2n}\}$ and $%
\alpha \neq 0.$
\end{lemma}

It follows that the set of characteristic vectors of $\frak{g}_{\Bbb{C}}$ is
the open set $\frak{g}_{\Bbb{C}}-\Bbb{C}\{X_{2},..,X_{2n}\}.$

\begin{lemma}
Either $\frak{g}_{1}$ or $\sigma (\frak{g}_{1})$ contains a characteristic
vector of $\frak{g}_{\Bbb{C}}.$
\end{lemma}

Observe that otherwise we would have $\frak{g}_{1}\subset \Bbb{C}$\{$%
X_{2},..,X_{2n}\}$ and $\sigma (\frak{g}_{1})\subset \Bbb{C}$\{$%
X_{2},..,X_{2n}\}$, which contradicts the previous decomposition.

Thus $\frak{g}_{1}$ or $\sigma (\frak{g}_{1})$ is a $n-$dimensional complex
filiform Lie algebra. But if $Y\in $ $\frak{g}_{1}$ is a characteristic
vector of $\frak{g}_{1}$, then $\sigma (Y)$ is a characteristic vector of $%
\sigma (\frak{g}_{1})$ with the same characteristic sequence.\ Then  every $%
2n-$dimensional filiform Lie algebra appears as a direct vectorial sum of
two $n$-dimensional filiform Lie algebras. We shall prove that it is
impossible. More precisely we have

\begin{proposition}
Let $n\geq 3$. Then no $2n-$dimensional complex filiform Lie algebra is a
vectorial direct sum of two $n-$dimension filiform sub-algebras.
\end{proposition}

{\bf Proof}. Let $\frak{g}_{\Bbb{C}}$ be a filiform Lie algebra of dimension $2n$
such that $\frak{g}_{\Bbb{C}}=\frak{g}_{1}\oplus \frak{g}_{2}.\;$From the
previous lemma, one of them for example $\frak{g}_{1}$, contains a
characteristic vector. Let be $X_{1}$ this vector an $%
\{X_{1},X_{2},..,X_{2n}\}$ the corresponding basis. This implies that $\frak{%
g}_{1}\cap \Bbb{C}\{X_{2},..,X_{2n}\}=\frak{g}_{1}\cap \Bbb{C}%
\{X_{n+1},..,X_{2n}\}.$ But $\frak{g}_{2}$ cannot contain characteristic
vector of $\frak{g},$ if not $\frak{g}_{2}\cap \Bbb{C}\{X_{2},..,X_{2n}\}=%
\frak{g}_{2}\cap \Bbb{C}\{X_{n+1},..,X_{2n}\}$ and this is at variance with $%
\frak{g}_{\Bbb{C}}=\frak{g}_{1}\oplus \frak{g}_{2}$. Then $\frak{g}%
_{2}\subset \Bbb{C}\{X_{2},..,X_{2n}\}.$ But from the brackets in $(*)$ this
is impossible as soon as $n>2.$

The four dimensional case can be treated directly. Up to isomorphism,  there
exists only one 4 dimensional filiform Lie algebra, $L_{4}$, for which we
have proven the nonexistence of invariant complex structures.\ Note that,
respect to the basis used in $(*)$, this algebra admits the decomposition $%
L_{4}=\frak{g}_{\Bbb{C}}=\frak{g}_{1}\oplus \frak{g}_{2}$,where $\frak{g}_{1}
$ and $\frak{g}_{2}$ are the abelian subalgebras generated respectively by $%
\{X_{1},X_{4}\}$ and $\{X_{2},X_{3}\}$.

\subsection{Consequence.}

Let $J$ be an invariant complex structure on a Lie algebra $\frak{g.}$ Let
us note by $\mu $ the law (the bracket) of $\frak{g}$ and let us consider
the Chevalley cohomology of $\frak{g.}$ The coboundary operator is denoted
by $\delta _{\mu }.$

\begin{proposition}
We have 
\[
\delta _{\mu }J=\mu _{J}
\]
where $\mu _{J}$ is the law of Lie algebra, isomorphic to $\mu $, defined by 
\[
\mu _{J}(X,Y)=J^{-1}(\mu (J(X),J(Y)).
\]
\end{proposition}

In fact the Nijenhuis condition is written as: 
\[
\mu (JX,JY)=\mu (X,Y)+J\mu (JX,Y)+J\mu (X,J(Y)) 
\]
Then 
\[
J^{-1}\mu (JX,JY)=J^{-1}\mu (X,Y)+\mu (JX,Y)+\mu (X,J(Y)) 
\]
that is, as $J^{2}=-Id$%
\begin{eqnarray*}
J^{-1}\mu (JX,JY) &=&-J\mu (X,Y)+\mu (JX,Y)+\mu (X,J(Y)) \\
&=&\delta _{\mu }J(X,Y)
\end{eqnarray*}

\begin{corollary}
If $\frak{g}$ is a filiform Lie algebra,then there does not exist
2-coboundaries for the Chevalley cohomology such that 
\[
\delta _{\mu }(J)=\mu _{J}
\]
where $J^{2}=-Id$
\end{corollary}

\[
\]

\end{document}